\documentclass[12pt]{amsart}
\usepackage{amsmath,amssymb,amsthm,amscd}
\numberwithin{equation}{section}

\def\SS{\Bbb S}

\def\H{{\Bbb H}}

\def\R{{\mathfrak R}}

\newtheorem{prop}{Proposition}[section]
\newtheorem{theo}[prop]{Theorem}
\newtheorem{lemm}[prop]{Lemma}
\newtheorem{coro}[prop]{Corollary}

\newtheorem{defi}[prop]{Definition}

\def\begeq{\begin{equation}}
\def\endeq{\end{equation}}
\def\p{\partial}

\def\lf{\left}
\def\ri{\right}
\def\e{\epsilon}
\def\ol{\overline}
\def\R{\Bbb R}

\begin{document}
\title{Asymptotically hyperbolic metrics on the unit ball with horizons}
\author{Yuguang Shi$^1$}
\address{Key Laboratory of Pure and Applied mathematics, School of Mathematics Science, Peking University,
Beijing, 100871, P.R. China.} \email{ygshi@math.pku.edu.cn}

\author{Luen-Fai Tam$^2$}

\thanks{$^1$ Research partially supported by NSF grant of China.
\\$^2$Research
partially supported by Earmarked Grant of Hong Kong \#CUHK403005}

\address{Department of Mathematics, The Chinese University of Hong Kong,
Shatin, Hong Kong, China.} \email{lftam@math.cuhk.edu.hk}

\date{May 2006}

\begin{abstract} In this paper, we construct a family of
asymptotically hyperbolic manifolds with horizons and with scalar
curvature equal to $-6$. The manifolds we constructed can be
arbitrary close to anti-de Sitter-Schwarzschild manifolds at
infinity. Hence, the mass of our manifolds can be very large or very
small. The main arguments we used in this paper is gluing methods
which was used in \cite{M}.
\end{abstract}

\maketitle \markboth{Yuguang Shi and Luen-Fai Tam} {On the
construction of asymptotically hyperbolic manifolds with horizons}

\section{Introduction}
In the past few years, there are many works  on the construction
of asymptotically flat (AF) and scalar flat manifolds which
contain minimal spheres. See \cite{BM}, \cite{M},   \cite{Co}, and
\cite{ST} for examples, and for existence of many blackholes,
please see \cite{CD}. From the point of view of general
relativity, these are examples of globally
 regular and asymptotically flat initial data for the Einstein vacuum equations
 containing a trapped surfaces.
 According to \cite{MSY},  if the topology of an AF manifold
 is nontrivial, then this  manifold always contains an outer most minimal
 sphere. The examples in \cite{BM},\cite{Co},\cite{M}, \cite{ST}, have the interesting property that
 the manifolds in those examples are all diffeomorphic to
$\R^3$.

 Another natural class of manifolds that are of interest in general
 relativity consists of asymptotically hyperbolic (AH) manifolds
 ( see Definition 1.1).
 Such manifolds arise when considering solutions to the Einstein
 fields equations with a negative cosmological constant, or when
 considering ``hyperboloidal hypersurfaces" in space-times which are
 asymptotically flat in isotropic directions. Therefore, it seems to
 be interesting to find    $3$-dimensional AH manifolds with
 $R=-6$ with trivial topology which  contain  horizons. In the
 asymptotically hyperbolic context,  horizons refer    not only
 to boundaries of domains which are minimal  but
 also to
 boundaries satisfying $H=\pm 2$. Here $H$ is the mean curvature of
 the boundaries with respect to the outward unit normal vectors.
 More precisely, we are   interested in following:

 {\it To
  find an AH manifold which is diffeomorphic to an open
 $3$-ball in $ \mathbf{R}^3$ with scalar curvature $R=-6$ which contains
  spheres with $H=0$ or $\pm 2$}.

 In  the AF context, in \cite{M},
 Miao constructs an AF and scalar flat manifold with topology
 $\mathbb{R}^3$ and containing a horizon (see also \cite{BM}).
The main arguments in \cite{M} is to glue  $\mathbb{S}^3$ with
Schwarzchild manifold and then conformally deform the metric to a
scalar flat  AF metric so that it still contains a minimal sphere.
We will use   similar methods
 to study our problem. More precisely, we will glue the
 anti-de Sitter-Schwarzchild space (see Section 1 for   details)
 with   part of the unit ball and obtain a complete metric with
 scalar curvature $R\ge -6$, which contains topological spheres with $H=0$
 or $\pm 2$, so that the metric is conformal to the hyperbolic metric
 on the unit ball in $\R^3$. Moreover, outside a compact set, the manifold
  is
 part of the anti-de Sitter-Schwarzchild space.
  Then we will    deform the metric to obtain
  an  AH with $R=-6$ which contains spheres with mean curvature 0,
  $\pm2$.
  We can show that the mass of our
manifolds (in the sense of \cite{Wang})can be close enough to that
of the anti-de Sitter-Schwarzchild space provided that the
perturbation is small enough. Hence the mass can be very large or
very small.

The   outline of the paper is as follows. In Section 1 we discuss
some basic facts of   anti-de Sitter-Schwarzchild space. Most of
them are well known, but we cannot find the details in literatures.
In Section 2, we will construct AH metric on the unit ball in
$\mathbb{R}^3$ which contains horizons. The metrics are rotationally
symmetric and are conformal to the hyperbolic metric, with scalar
curvature $R\ge-6$ so that  $R=-6$ near infinity. In Section 3, we
will do the deformation  to obtain new AH metrics on the ball with
scalar curvature equal to $-6$ which contain horizons. We also
discuss  the mass  of these AH metrics in this section.

We would like to thank Lars Andersson for useful discussion which
motivates this work.

\section{The anti-de Sitter-Schwarzschild metric}

In order to construct metrics which are asymptotically
hyperbolic (AH) and contains horizons, we will make use of
the anti-de Sitter-Schwarzschild metric. Therefore in this
section, we will discuss this metric in details. Let us
first recall the definition of asymptotically hyperbolic
manifolds and its mass. We will use the definitions in
\cite{Wang}, see also \cite{CH}. We are only interested in
the case that the manifold has dimension three.
\begin{defi} A complete noncompact Riemannian manifold
$(X^3,g)$ is said to be {\it asymptotically hyperbolic} if there is
a compact manifold $(\ol X,\ol g)$ with boundary $\p X$ and a smooth
function $t$ on $\ol X$ such that the following are true:
 \begin{enumerate}
   \item[(i)] $X=\ol X\setminus \p X$.
   \item [(ii)]$t=0$ on $\p X$, and $t>0$ on $X$.
   \item[(iii)] $\ol g=t^2g$ extends to be $C^3$ up to the
   boundary.
   \item[(iv)] $|d t|_{\ol g}=1$ at $\p X$.
   \item [(v)]Each component $\Sigma$ of $\p X$ is the standard two
   sphere $(\mathbb{S}^2,g_0)$ and there is a collar
   neighborhood of $\Sigma$ where
   $$
   g=\sinh^{-2}t(dt^2+g_t)
   $$
   with
   $$
   g_t=g_0+\frac{t^3}{3}h+O(t^{4})
   $$
   where $h$ is a $C^2$  symmetric two tensor on
   $\mathbb{S}^2$.
 \end{enumerate}

\end{defi}
With the above notation, let $(X,g)$ be an AH manifold with scalar
curvature $R\ge -6$, then the {\it mass} of an end of $X$
corresponding to a boundary component $\Sigma$ of $\p X$ is defined
as
$$
M=\frac1{16\pi}
\lf[\lf(\int_{\mathbb{S}^2}\text{trace}_{g_0}(h)dV_{g_0}\ri)^2
-\lf(\int_{\mathbb{S}^2}\text{trace}_{g_0}(h)(x)xdV_{g_0}\ri)^2
\ri]^\frac12
$$
where $x$ is the standard coordinates of a point on $\Bbb
S^2$ in $\Bbb R^3$. This is well-defined by \cite{Wang}.

Next we want to describe the anti-de Sitter-Schwarzschild
metric which is obtained by gluing two copies of manifolds
with boundary with metric $$ ds^2=\frac{dr^2}{1+r^2-\frac
Mr}+r^2d\sigma^2
$$
with $M>0$ defined on $(a(M),\infty)\times \SS^2$ where $a(M)>0$
is the unique root  of $1+r^2-\frac Mr=0$ and $d\sigma^2$ is the
standard metric on the standard sphere $\SS^2$. The construction
and some properties of the metric are well-known. But for the sake
of completeness and for reference later, we will give  details of
the metric and its properties.

Let
$$
h(r)=\int_r^\infty\frac{1}{\sqrt{t(t+t^3-M)}} dt
$$
for $r>a(M)$. Let $\rho(r)$ be the function defined by the
relation

$$
e^\rho=\frac{1+e^{-h}}{1-e^{-h}}
$$
that is
$$
\rho(r)=\log\lf(\coth \frac{h(r)}2 \ri).
$$
 Then $\rho:(a(M),\infty)\to ( \rho(M), \infty)$ is a
 smooth increasing function in $r$ with $\rho(M)=\rho(a(M))>0$.
 Let
 $\phi>0$ be the smooth function in $\rho$ on $( \rho(M),
 \infty)$ by
 \begin{equation}\label{phi}
\phi^4(\rho)=\frac{r^2(\rho)}{\sinh^2\rho}.
\end{equation}
We have
 \begin{equation}\label{1stmetric}
ds^2=\phi^4 (d\rho^2+\sinh^2\rho d\sigma^2).
\end{equation}
Here $d\rho^2+\sinh^2\rho d\sigma^2$ is the standard metric
on the hyperbolic space $\H^3$. Observe that $\phi$ is
continuous  up to $\rho(M)$ and is positive at $\rho(M)$.
Also,
\begin{equation}
  \phi_\rho =-\frac12 \phi^{-1}(\sinh h+rh_r\cosh h)\frac{\sinh
h}{h_r}
\end{equation}
and
\begin{equation}
\begin{split}
\phi_{\rho\rho}&=-\frac14\phi^{-3}\lf[(\sinh h+rh_r\cosh
h)\frac{\sinh h}{h_r}\ri]^2\\
&+\frac12\phi^{-1}\lf[3 \sinh h\cosh h-\frac{\sinh^2 h
h_{rr}}{h^2_r} +r(\sinh^2 h+\cosh^2h)h_r\ri]\frac{\sinh
h}{h_r}
\end{split}
\end{equation}
From these, it is easy to see that $\phi$ as a function of
$\rho$ is $C^2$ up to $\rho(M)$. The scalar curvature of
$ds^2$ is:
\begin{equation}\label{scalar}
    \begin{split}
R&=\phi^{-5}\lf(-6\phi-8\Delta_{\H^3}\phi\ri)\\
&=\phi^{-5}(-6\phi-8(\phi_{\rho\rho}+2\coth\rho\phi_\rho))\\
&=2\lf[ h_r^{-2} r^{-4} + r^{-2} +2h_{rr}h_r^{-3} r^{-3}
)
\ri]\\
&=-6
\end{split}
\end{equation}

Let us use the ball model for $\H^3$. Let
$$
b(M)=\frac{e^{\rho(M)}-1}{e^{\rho(M)}+1}.
$$
Then
\begin{equation}\label{2ndmetric}
ds^2=\frac{4\psi^4(\tau)}{(1-\tau^2)^2}(d\tau^2+\tau^2d\sigma^2),
\end{equation}
on the annulus $b(M)\le|x|<1$ in $\R^3$, where $\tau=|x|$
and
$$\psi(\tau)=\phi\lf(\log\lf(\frac{1+\tau}{1-\tau}\ri)\ri).$$

Now we use the transformation $x\to b^2(M)x/|x|^2$ to
transform the annulus $b(M)\le |x|< 1$ to $b^2(M)< |x|\le
b(M)$, such that $|x|=1$ is mapped to $|x|=b^2(M)$. Pull
back the metric $ds^2$ to $b^2(M)< |x|< b(M)$, we extend
the metric  $ds^2$ to  $b^2<|x|<1$ and is still denoted by
$ds^2$ such that the metric is of the form
$$
ds^2=f^4(\tau)ds^2_{\Bbb H^3},
$$
where
$$
f(\tau)=\left\{%
\begin{array}{ll}
    & \psi(\tau)\hbox{\qquad for $b\le |x|<1$;} \\
      &b\psi(\frac{b^2}\tau)\lf(\frac{1-\tau^2}{\tau^2-b^4}\ri)^\frac12
       \hbox{\ for $b^2< |x|\le b$.} \\
\end{array}%
\right.
$$
where $|x|=\tau$. It is easy to see that $f$ is continuous
at $b(M)$. We want to prove that $f_\tau$ matches at
$\tau=b(M)$. Suppose  $\tau\to b(M)_+$, then
\begin{equation}\label{1}
f_\tau=\psi_\tau\to \phi_\rho \rho_\tau=-\frac{\phi\cosh
\rho(M)}{\sinh
\rho(M)}\cdot\frac{1}{1-b^2(M)}=-\frac{\psi(1+b^2(M))}{2(1-b^2(M))}.
\end{equation}
Suppose $\tau\to b(M)-$, then
\begin{equation*}
    \begin{split}
f_\tau&\to b(M)\psi_\tau(b(M))(-\frac{b^2(M)}{\tau^2})\lf(\frac{1-\tau^2}{\tau^2-b^4(M)}\ri)^\frac12\\
&\quad+b(M)\psi(b(M))\frac12\lf(\frac{1-\tau^2}{\tau^2-b^4(M)}\ri)^{-\frac12}
\cdot\frac{(\tau^2-b^4(M))(-2\tau)-(1-\tau^2)\cdot2\tau}{(\tau^2-b^4(M))^2}\\
&=-\psi_\tau -\frac{\psi(1+b^2(M))}{b(M)(1-b^2(M))}\\
&=\psi_\tau
\end{split}
\end{equation*}
by (\ref{1}). Hence $f$ is $C^1$ near $\tau=b(M)$. Since
$ds^2$ has scalar curvature $-6$ on $b^2(M)<|x|<b(M)$ and
$b(M)<|x|<1$, $f$ is a $C^1$ weak solution of
\begin{equation}
\Delta_{\H^3}f-\frac34f(f^4-1)=0.
\end{equation}
By regularity, $f$ must be smooth. Hence we obtain a smooth
metric
$$
g_{_{AdS-Sch,M}}=\phi_M^4(d\rho^2+\sinh^2\rho d\sigma^2)
$$
defined on $\H^3\setminus B(\rho_0(M))$ where
$\phi_M(\rho)=f(\tau(\rho))$ and
$$\rho_0(M)=\log\lf(\frac{1+b^2(M)}{1-b^2(M)}\ri).
$$
Here $B(\rho_0(M))$ is the geodesic ball with center at
$\rho=0$. The metric $g_{_{\text{AdS-Sch},M}}$ is called
the anti-de Sitter-Schwarzschild metric (with mass $M$).

With the above notation, we have:
\begin{prop}\label{AdsSch1} For each $M>0$, the anti-de  Sitter-Schwarzschild metric
$g_{_{\text{AdS-Sch},M}}=\phi^4_M(d\rho^2+\sinh^2\rho
d\sigma^2)$ is complete and is defined on $\H^3\setminus
B(\rho_0(M))$.   Moreover:
\begin{enumerate}
  \item [(i)] $\phi>1$, $\lim_{\rho\to\infty}\phi=1$  and
  $\lim_{\rho\to\rho_0(M)}\phi=\infty$.
  \item [(ii)] The manifold is asymptotically hyperbolic with constant
  scalar curvature -6.
  \item[(iii)] Denote $\phi_M$ by $\phi$, then $\phi_\rho<0$ and
  \begin{equation}\label{metric1}
  (\sinh^2\rho \cdot\phi_\rho)_\rho=\frac34\sinh^2\rho\cdot\phi(\phi^4-1).
  \end{equation}

  \item [(iv)] There exist unique $\rho_2>\rho_1>\rho_2'>\rho_0(M)$
  such that the level surface of $\rho=\rho_1$ is minimal, the mean curvature of
  $\rho=\rho_2$ is 2 and $\rho=\rho_2'$ is -2 with respect to the unit
   normal
  in the direction $\frac{\p}{\p \rho}$.
\end{enumerate}
\end{prop}
\begin{proof} (i) The results are immediate from the
definition of $\phi_M$.

(ii) This follows from (\ref{scalar}) and \cite{Wang}.

(iii) $\phi $ satisfies (\ref{metric1}) because the scalar
curvature is -6. From the equation, we have
$\sinh^2\rho\phi_\rho$ is strictly increasing. Suppose
$\phi_\rho\ge 0$ for some $\rho^*$, then $\phi_\rho>0$ for
all $\rho>\rho^*$. Since $\phi>1$, this contradicts (i).
 (iv) Denote $\phi_M$ simply by $\phi$. The mean
curvature of the level surface $\rho=$constant for
$\rho>\rho(M)$ is
\begin{equation}\label{meancurvature}
\begin{split}
 H&=\frac1{\phi^2}\lf(\frac{2\cosh\rho}{\sinh
\rho}+\frac4\phi\phi_\rho\ri)\\
&=\frac1{\phi^2}\lf(2\cosh h-2\phi^{-2}(\sinh h+rh_r\cosh
h)\frac{\sinh h}{h_r} \ri)\\
&=-2\phi^{-4}\sinh^2 h h_r^{-1}\\
&=-2r^{-2}h_r^{-1}\\
&=2\lf(1+r^{-2}-Mr^{-3}\ri)^\frac12.
\end{split}
\end{equation}
From this the results follow.
\end{proof}

Next we will discuss the behaviors of the metrics
$g_{_{\text{AdS-Sch},M}}$ as $M$ changes. Before we do
this, we need the following lemma which may be well-known:
\begin{lemm}\label{maximum} Let $(N,g)$ be a complete noncompact Riemannian
 manifold. Suppose $u_1\geq 1$ and $u_2\ge 1$ are   such that
\begin{equation}\label{max}
    \Delta u_1 + \frac34 u_1(1-u_1^4)=\Delta u_2 + \frac34 u_2(1-u_2^4)
\end{equation}
 on
$N \setminus B(\rho^*)$ where $\Delta$ is Laplacian of $N$ and
$B(\rho^*)$ is the geodesic ball of radius $\rho^*$ with center at a
fixed point. Suppose $u_1\ge u_2$ at $\p B(\rho^*)$ and suppose
$\lim_{x\to\infty}(u_1(x)-u_2(x))=0$, then $u_1\ge u_2$ in $N
\setminus B(\rho^*)$. If in addition that $N= \H^3$, then
$$
|u_1-u_2|(x)\le C\lf(\sup_{\p B(\rho^*)}|u_1-u_2|\ri) e^{-3\rho(x)}
$$
outside $B(\rho^*)$ where $\rho$ is the distance function from a
fixed point and $C$ is a constant depending only on $\rho^*$. In
case $\rho^*=0$, then $u_1=u_2$.
\end{lemm}
\begin{proof} Let us prove the last statement and the first assertion can be proved similarly.
 Let $\eta=u_1-u_2$, then
$$
\Delta \eta=3\eta+\frac34\eta (-5+G) $$
 where $G=u_1^4+u_1^3u_2+u_1^2u_2^2+u_1u_2^3+u_2^4\ge5$.  Let $\xi(\rho)=
e^{-2\rho}\sinh^{-1}\rho$, then it is easy to check that
$$\Delta \xi = 3 \xi.$$
Let $$A=\frac{ \sup_{\p B(\rho^*)}|u_1-u_2| }{\xi(\rho^*)}.$$ Then
by maximum principle, we can conclude that
$$
\eta\le A\xi
$$
outside $B(\rho^*)$. Similarly, one can prove that
$-\eta\le A\xi$. From this the second part of the lemma is
proved.
\end{proof}

 Now we are ready to  discuss the behaviors of the metrics
$g_{_{\text{AdS-Sch},M}}$. Let $a(M)$, $\rho(M)$ and $b(M)$
be as before. The metric $g_{_{\text{AdS-Sch},M}}$  is of
the form $\phi_{_{M}}^4ds^2_{\H^3}$ which is defined and is
complete on $b^2(M)<\tau<1$ in the ball model of $\H^3$.
\begin{prop}\label{massmetrics} With the above notation, we have
the following:

\begin{itemize}
    \item [(i)] $a(M)$, $\rho(M)$, $b(M)$ are continuous monotonic increasing
     functions of $M$.
     \item[(ii)]
      \begin{equation}\label{massmetrics1}
\lim_{M\to 0}\frac{a(M)}{M}=1;\  \lim_{M\to
0}\rho(M)=\lim_{M\to
    0}b(M)=0.
\end{equation}
\begin{equation}\label{massmetrics2}
\lim_{M\to\infty}\frac{a}{M^\frac13}=1;\
     \lim_{M\to \infty}\rho(M)=\infty;\
     \lim_{M\to\infty}b=1.
\end{equation}
    \item [(iii)]   $M_1>M_2>0$ if and only if $\phi_{M_1}>\phi_{M_2}$ on
    $b^2(M_1)<\tau<1$.
    \item [(iv)] For each fixed $0<\tau<1$,  $\phi_M(\tau)$ is a
    continuous function of $M$ whenever it is defined  and
    $\lim_{M\to0}\phi_M(\tau)=1$.
\end{itemize}
\end{prop}
\begin{proof} (i) It is easy to see that $M_1>M_2$ implies
$a(M_1)>a(M_2)$ and $a(M)$ is continuous in $M$. Next we
want to prove that $h(a(M_1))<h(a(M_2))$. Given $M>0$, let
us denote $a=a(M)$ for simplicity. Then
$t+t^3-M=(t-a)(t^2+at+1+a^2)$ by direct computation. Hence
\begin{equation}\label{massmetrics3}
    \begin{split}
h(a(M))&=\int_a^\infty \frac{1}{\sqrt{t(t-a)(t^2+at+1+a^2)}}dt\\
&=\int_0^\infty \frac{1}{\sqrt{t(t+a) (t^2+3at+1+3a^2)}}dt.
\end{split}
\end{equation}
From this it is easy to see that $h(a(M_1))<h(a(M_2))$ if
$M_1>M_2$. Hence $\rho(M_1)>\rho(M_2)$ and $b(M_1)>b(M_2)$.
From (\ref{massmetrics3}) it is easy to see that $h(a(M))$
is a continuous function of $a(M)$ and hence is continuous
in $M$. So $\rho(M)$ and $b(M)$ are continuous in $M$. This
proves (i).

(ii) Given $M>0$ denote $a(M)$ simply by $a$. Then
$a+a^3-M=0$ and so $a<M$ and   $a\to0$ as $M\to0$.
$$
M>a=M-a^3>M-M^3.
$$
  From this, we can conclude that
$a/M\to 1$ as $M\to 0$. By (\ref{massmetrics3}), we have
$h(a(M))\to\infty$ if $M\to0$ (and so $a(M)\to 0$). Hence
(\ref{massmetrics1}) is true.

It is easy to see that if $M\to\infty$ then
$a=a(M)\to\infty$ and $a^3<M$. On the other hand, let
$1>\delta>0$, then $a^3\ge M-\delta a^3$ provided $M$ is
large enough. From this, (\ref{massmetrics2}) follows.

(iii) Suppose $M_1>M_2$, then $b^2(M_1)>b^2(M_2)$. Hence
$\phi_{M_2}$ is bounded on $ \tau=b^2(M_1)$ and $\phi_{M_1}=\infty$
at $ \tau=b^2(M_1)$. Since both $\phi_{M_1}>1$ and $\phi_{M_2}>1$
satisfies the equation: $\Delta u + \frac34 u (1-u ^4)=0$ outside
the geodesic ball in $\H^3$ corresponding to $|x|<b^2(M_1)$ in the
ball model, and since $\phi_{M_1}, \phi_{M_2}\to 1$ as $\tau\to1$,
(iii) follows from Lemma \ref{maximum}.

(iv) Let  $\tau$ be fixed. For any $M_0$ such that $\tau>b^2(M_0)$
then $\tau>b^2(M)$ provided $M$ is close enough to $M_0$. By the
construction of $\phi_M$, it is sufficient to prove that case that
$\tau>b(M_0)$. By the construction, it is sufficient to prove the
following: if $\rho>0$ is fixed such that $\rho>\rho(M_0)$, then
$\phi_M(\rho)\to \phi_{M_0}(\rho)$ as $M\to M_0$. Now by
(\ref{phi}),
$$
\phi_M^4(\rho)=\frac{r^2}{\sinh^2\rho}.
$$
where $r$ and $\rho$ is related by
$$
\sinh h(r)=\frac1{\sinh \rho}
$$
with
$$
h(r)=\int_r^\infty\frac{1}{\sqrt{t(t+t^3-M)}} dt.
$$
From these it is easy to see the result follows.

To prove the second assertion in (iv), it is sufficient to
prove that for fixed $\rho$, $r\sinh h(r)\to 1$ as $M\to0$,
where $r$ and $h(r)$ are given by $\sinh h=1/\sinh\rho$ and
$$
h(r)=\int_r^\infty \frac{dt}{\sqrt{t(t+t^3-M)}}.
$$
Then as $M\to 0$, $r\to r_0$ such that $\sinh
h(r_0)=1/\sinh\rho$ and
$$
h(r_0)=\int_{r_0}^\infty \frac{dt}{\sqrt{t(t+t^3)}}.
$$
Hence $\sinh h(r_0)=1/r_0$ and so $r_0=\sinh \rho$. From
(\ref{phi}), the result follows.
\end{proof}

As an application, we have the following uniqueness result:
\begin{coro}\label{uniqueness} Suppose $g=\phi^4ds^2_{\H^3}$ is a conformal metric
defined on $\H^3\setminus B_\rho$ for some $\rho>0$ such
that the scalar curvature is -6. Suppose $\lim_{x\to
\infty}\phi(x)= 1$, $\phi>1$ and $\phi=$constant on $\p
B_\rho$. Then $g=g_{_{\text{AdS-Sch},M}}$ on $\H^3\setminus
B_\rho$ for some $M>0$.
\end{coro}
\begin{proof} The corollary follows from Lemma \ref{maximum}, Propositions
\ref{AdsSch1} and   \ref{massmetrics}.
\end{proof}

\section{ Conformal AH metrics on the unit ball}

In this section, we will construct asymptotically
hyperbolic (AH) metrics on the unit ball in $\R^3$  which
contains horizons and which is conformal to the hyperbolic
metric. Moreover, the  scalar curvature $R$ satisfies $R\ge
-6$ and the manifold is a part of the anti-de
Sitter-Schwarzschild near infinity.

Let $M>0$ and let
 $g_{_{\text{AdS-Sch},M}}=\phi^4_Mds^2_{\H^3}$  be the
anti-de Sitter -Schwarzschild metric defined in Proposition
\ref{AdsSch1}. Let $\rho_2>\rho_1>\rho_2'>\rho_0(M)>0$ be as in
the proposition.   First we want to construct a $C^{2,1}$ metric
with the  properties mentioned above such that it is anti-de
Sitter-Schwarzschild outside $B(\tau_2)$ for some
$\rho_2'>\tau_2>\rho_0(M)$, where  $B(\tau_2)$ is the geodesic
ball with center at $\rho=0$ of the hyperbolic space with metric
of the form $d\rho^2+\sinh^2\rho d\sigma^2$. Let us denote
$\phi_M$ simply by $\phi$. Note that if $f^4ds^2_{\H^3}$ is a
conformal metric such that $f$ depends only on $\rho$, then the
scalar curvature is given by
\begin{equation}\label{scalar1}
   R=f^{-5}(-6f-\Delta_{\H^3}f)=f^{-5}\lf[-6f
   -8\lf(f_{\rho\rho}+2\coth\rho\cdot f_\rho\ri)\ri].
\end{equation}

\begin{lemm}
\label{AdsSch2} With the above notations, there exist
$\rho_0(M)<\tau_1<\tau_2<\rho_2'$ and   a $C^{2,1}$ function
$\psi(\rho)$ on $[0,\infty)$ such that $\psi(\rho)>1$,
$\psi(\rho)=$constant on $[0,\tau_1]$ and such that
$\psi(\rho)=\phi(\rho)$ on $[\tau_2,\infty)$. Moreover, the metric
$\psi^4ds^2_{\H^3}=\psi^4(d\rho^2+\sinh^2\rho d\sigma^2)$ has
scalar curvature $R>-6$ on $B({\tau_2})$.
\end{lemm}
\begin{proof} For any
$\rho_0(M)<\tau_1<\tau_2<\rho_2'$, let $
\xi(\rho)=(\rho-\tau_1)^2(a\rho+b)$ where $a$ and $b$ are chosen so
that
\begin{equation}\label{metric2}
\left\{\begin{array}{ll}
\xi(\tau_1)&=0\\
\xi_\rho(\tau_1)&=0; \\
\xi(\tau_2)&=A=\sinh^2\tau_2\cdot\phi_\rho(\tau_2);\\
\xi_\rho(\tau_2)&=B=\frac34\sinh^2\tau_2\cdot\phi(\tau_2)(\phi^4(\tau_2)-1).
\end{array}
\right.
\end{equation}
Then $$a= (\tau_2-\tau_1)^{-2}\lf[B-2A(\tau_2-\tau_1)^{-1}\ri]$$
and
$$b=(\tau_2-\tau_1)^{-2}\lf[A-
\tau_2 B+2A\tau_2(\tau_2-\tau_1)^{-1}\ri].$$ Since $A<0$ and
$B>0$, we have $a>0$ and $b<0$. Since $\xi(\tau_2)=A<0$, so
$a\tau_2+b<0$. Since $a>0$, we have $a\rho+b<0$ for all
$\rho<\tau_2$. In particular,
$$
\xi\le0
$$
on $[\tau_1,\tau_2]$.

Define $\psi$ as follows

\begin{equation}\label{psi}
  \psi(\rho)=  \left\{%
\begin{array}{ll}
 \phi(\tau_2)-\int_{\tau_1}^{\tau_2}\frac{\xi(t)}{\sinh^2t}dt,
&\hbox{$0\le \rho<\tau_1$ ;} \\
    \phi(\tau_2)-\int_\rho^{\tau_2}\frac{\xi(t)}{\sinh^2t}dt,
    &\hbox{ $ \tau_1\le \rho\le\tau_2 $ ;} \\
  \phi(\rho),  &\hbox{$   \tau_1<\rho<\infty)$.} \\
\end{array}%
\right.
\end{equation}
Since in $[\tau_1,\tau_2]$, $\sinh^2\rho\cdot
\psi_\rho(\rho)=\xi(\rho)$, and since the metric
$\phi^4ds^2_{\H^3}$ has constant curvature $-6$, by Proposition
\ref{AdsSch1} and the definition of $\xi$, one can see that $\psi$
is $C^{2,1}$.

We want to compute the scalar curvature of $\psi^4ds^2_{\H^3}$.
Since $\xi<0$, $\psi>1$, then the scalar curvature on
$B({\tau_1})$ is larger than $-6$ because $\psi>1$ and is constant
there. Outside $B({\tau_2})$, $\psi=\phi$ and the scalar curvature
is -6. In $[\tau_1,\tau_2]$,
 $$
 \xi_\rho(\rho)=3a\rho^2+2(b-2a\tau_1)\rho+(a\tau_1^2-2b\tau_1).
 $$
Hence for  $\rho\in[\tau_1,\tau_2]$,
\begin{equation*}
    \begin{split}
     \xi_\rho(\tau_2)- \xi_\rho(\rho)
     &=3a(\tau_2^2-\rho^2)+2(b-2a\tau_1)(\tau_2-\rho)\\
     &=a(\tau_2-\rho)\lf[3(\tau_2+\rho)+\frac{2b}{a}-4\tau_1\ri]\\
     &\ge
     a(\tau_2-\rho)\lf[\tau_2-\tau_1+\frac{2A}{B-2A(\tau_2-\tau_1)^{-1}}\ri]\\
     &=a(\tau_2-\rho)\frac{B(\tau_2-\tau_1)^2}{B(\tau_2-\tau_1)-2A}\\
&\ge0
     \end{split}
\end{equation*}
and is positive if $\rho<\tau_2$, because $a>0$, $B>0$ and $A<0$.

  So  in  $B({\tau_2})\setminus
B({\tau_1})$, we have
\begin{equation*}
   \begin{split}
      \lf[\sinh^2\rho
\psi_\rho(\rho)\ri]_\rho &=\xi_\rho(\rho) \\
         & < \xi_\rho(\tau_2)\\
         &=\frac34\phi(\tau_2)(\phi^4(\tau_2)-1)\sinh^2\tau_2\\
         &\le\frac34\phi(\rho)(\phi^4(\rho)-1)\sinh^2\rho\\
         &\le\frac34\psi(\rho)(\psi^4(\rho)-1)\sinh^2\rho
     \end{split}
\end{equation*}
and the scalar curvature is larger than $-6$ by (\ref{scalar1}),
provided that $\phi(\rho)(\phi^4(\rho)-1)\sinh^2\rho$ is
decreasing on $[\tau_1,\tau_2]$. Here we have used the fact that
$\psi\ge \phi>1$ on $[\tau_1,\tau_2]$.

Now
\begin{equation*}
    \begin{split}
       \lf[\log\lf(\phi(\rho)(\phi^4(\rho)-1)\sinh^2\rho\ri)\ri]_\rho  &
        =\frac{2\cosh \rho}{\sinh\rho}+\frac{5\phi^4-1}{\phi^4-1}
        \lf(\log\phi\ri)_\rho \\
         &=\frac{2\cosh \rho}{\sinh\rho}+\lf(\log\phi\ri)_\rho
         +\frac{4\phi^4}{\phi^4-1}
        \lf(\log\phi\ri)_\rho.
     \end{split}
\end{equation*}
Since $\phi\to\infty$ as $\rho\to\rho_0(M)_+$ by Proposition
\ref{AdsSch1}, there exists $\rho_0(M)<\tau_2<\rho_2'$ such that
\begin{equation}\label{metric3}
\lf[\log\lf(\phi(\rho)(\phi^4(\rho)-1)\sinh^2\rho\ri)\ri]_\rho<0
\end{equation}
at $\tau_2$. Hence one can choose $\alpha<\tau_1<\tau_2<\rho_0$
such that (\ref{metric3}) is true in $[\tau_1,\tau_2]$. This
completes the proof of the lemma.
\end{proof}
Next we want to modify $\psi$ in the lemma so that it is
smooth.
 Using the same notation as in Lemma \ref{AdsSch2}. Let
 $$
 f=\Delta_{\H^3}
\psi+\frac34\psi(1-\psi^4).
$$
  Then $f$ is Lipschitz and
$f=f(\rho)<0$ on $[0,\tau_2)$ and $f=0$ on
$[\tau_2,\infty)$. For any $\e>0$ let $0\le \chi_\e\le 1$
be a cutoff function on $[0,\infty)$ such that $\chi_\e=1$
on $[0, \tau_2-\e)$ and $\chi_\e=0$ on
$[\tau_2-\frac12\e,\infty)$. Define $f_\e=f\chi_\e$. Then
$f_\e\ge f$ and $f_\e-f\le C(\e)$ where $C(\e)$ is a
function of $\e$ with $\lim_{\e\to0}C(\e)=0.$ Note that
$f_\e$ is smooth and $f_\e=f=0$ on $[\tau_2,\infty)$.

We want to prove the following:

\begin{theo}\label{AdsSch3} For any $\tau_2>\e>0$, there is a
unique $\phi_\e$ which depends only on $\rho$ such that
\begin{itemize}
    \item[(i)] \begin{equation}\label{metric4}
\Delta_{\H^3} \phi_\e+\frac34\phi_\e(1-\phi_\e^4)=f_\e,
\end{equation}
and hence the scalar curvature of the metric
$g_\e=\phi_\e^4ds^2_{\H^3}$ is not less  than -6 in $B({\tau_2})$
and is -6 outside $B({\tau_2})$.
    \item [(ii)]$$
\psi>\phi_\e>1
$$
and $\lim_{\rho\to\infty}\phi_\e=1$
    \item[(iii)] $\psi(x)-\phi_\e(x)\le C(\e)e^{-3\rho(x)}$ in
    $\H^3\setminus B(\tau_2)$, where $C(\e)\to0$ as $\e\to0$.
    \item[(iv)]
    $g_\e^4ds^2_{\H^3}=g_{_{\text{Ads-Sch},M_\e}}$
    for some $M_\e>0$ on $\H^3\setminus B(\tau_2)$.
     In particular, $g_\epsilon$ is AH. Moreover, if $\e>0$ is small enough, then
      $\tau_2<\rho_{2,\epsilon}'$
     where $\rho=\rho_{2,\e}'$ is the surface with constant mean curvature
      $-2$ in the metric $g_{_{\text{Ads-Sch},M_\e}}$.
    \item[(v)] Let $M_\e$ be as in (iv), then $M-M_\e\le
    C(\e)$, where $C(\e)\to0$ as $\e\to0$.
\end{itemize}
\end{theo}

\begin{proof} The existence part follows from \cite{AM}. In fact, let $\psi$ as in  Lemma
\ref{AdsSch2}. By the definitions of $f$ and $f_\e$, we have
$$
\Delta_{\H^3}\psi+\frac34\psi(1-\psi^4)=f\le f_\e
$$
and $\psi>1$, $\psi\in C^{2,1}_{loc}(M)$. Since $f_\e\le0$, we
have
$$
\Delta_{\H^3}\psi_0+\frac34\psi_0(1-\psi^4_0)=0\ge f_\e
$$
where $\psi_0=1$ is the constant function. By \cite{Sa}, for any
integer $k\ge 1$, we can find a unique solution $\psi_k$
\begin{equation}\label{metric6}
\left\{
\begin{array}{ll}
 \Delta_{\H^3} \psi_k +\frac34 \psi_k(1-{\psi_k}^4)=f_\e, & \hbox{ in $B(k)$;} \\
   \psi_k |_{\p B(k)}=\psi |_{\p B(k)},  \\
\end{array}
\right.
\end{equation}
with $1\le \psi_k\le \psi$. Hence one can choose a subsequence of
$\psi_k$ which converges uniformly on compact subsets of $\H^3$
together with its derivatives to a solution $\phi_\e$ of
(\ref{metric4}). Moreover, $1\le \phi_\e\le\psi$ and hence
$1<\phi_\e<\psi$ by the strong maximum principle. Moreover, since
$f_\e$ is a function of  $\rho$, $\phi_\e$ is also a function of
$\rho$ by Lemma \ref{maximum}. This proves (i) and (ii).

For $k>\tau_2$, let  $\eta=\psi_k -\psi$. Then
$$
\Delta_{\H^3} \eta+\frac34\eta(1-G)=f_\e-f
$$
 where $G=\psi^4 +\psi^3 \psi_k + \cdots+ \psi_k ^4 >5$. Multiply
 both sides by $\eta$  and integrating by parts, we
get

$$\int_{B(k)}|\nabla \eta|^2 + \frac34 \int_{B(k)} G \eta^2 -\frac34 \int_{B(k)} |\eta|^2
=- \int_{B(k)}(f_\e -f)\eta=- \int_{B(\rho_2)}(f_\e -f)\eta.$$ Hence
there exists a function $C(\e)$ such that
$$
\int_{B(k)}|\eta|^2\le C(\e)
$$
here and below, $C(\e)$ denotes a function of $\e$ such that
$\lim_{\e\to0}C(\e)=0$

 Hence we have
$$
\int_M|\phi_\e-\psi|^2\le C(\e).
$$
By mean value inequality \cite{GT}, we conclude that for any
$\rho>\tau_2$,
$$
\sup_{B(\rho) \setminus B(\tau_2)}|\phi_\e-\psi|\le C(\e).
$$
Since both $\psi$ and $\phi_\e$ satisfy
$$
\Delta_{\H^3} u+\frac34(1-u^4)=0
$$
on $\H^3\setminus B({\tau_2})$, by Lemma \ref{maximum}, we
conclude that
$$
\sup_{\H^3\setminus B({\tau_2})} \psi(x)-\phi_\e (x)\le
C(\e)e^{-3\rho(x)}.
$$
This proves (iii).

The first part of (iv) follows from Corollary \ref{uniqueness}.
(v) follows from (ii), (iii) and   Proposition \ref{massmetrics}.
The second part of (iv) follows from (v) and Proposition
\ref{massmetrics}

\end{proof}

\section{AH metrics with $R=-6$ on the unit ball  with horizons}

Using the metrics constructed in \S2, we will   construct AH
metrics
 on the unit ball with $R=-6$ which contains a minimal sphere and
 spheres with constant mean curvature $\pm2$. More precisely, we
 have:
 \begin{theo} Let $\mathbf{D}$ be the unit ball in $\R^3$. For any
 $M>0$ and $\delta>0$, there is a smooth complete metric $g$ on
 $\mathbf{D}$ with constant scalar curvature $-6$ such that the
 following are true
  \begin{enumerate}\label{AdsSch4}
    \item [(i)] $(\mathbf{D},g)$ is asymptotically hyperbolic with mass $M_g$ satisfying
    $|M_g-M|<\delta$.
    \item [(ii)] There exist surfaces $S_1$, $S_2$, and $S_3$
    which are topological spheres with constant mean curvature $-2, 0, 2$ respectively
    such that $S_1$ is in the
    interior of $S_2$ and $S_2$ is in the interior of $S_3$.
    \item[(iii)] Outside a compact set the metric $g$ is conformal to the
    standard hyperbolic metric of $\mathbf{D}$.
 \end{enumerate}
 \end{theo}
\begin{proof}
Let $M>0$ and $\e>0$ be given, let   $g_1=g_\e=\phi_\e
^4ds^2_{\H^3}$ be the metric constructed in Theorem \ref{AdsSch3}.
The scalar curvature $R_1$ of $g_1$ is $-6$ outside the geodesic
ball $B(\tau_2)$ with respect to $ds^2_{\H^3}$, and $R_1\ge -6$.
Let   $0\le \xi\le 1$ be a smooth function which is positive on
$B(\tau_2)$, zero outside $B(\tau_2)$. For $\delta>0$, let
$f_\delta=-6-\delta\xi$. Then $f_\delta<R_1$ in $B(\tau_2)$ and
$f_\delta=R_1$ outside $B(\tau_1)$. By the result of Lohkamp
\cite[Theorem 1]{Lo}, there is a metric $g_2$ such that the scalar
curvature $R_2$ of $g_2$ satisfies $f_\delta-\delta\le R_2\le
f_\delta$ on $B(\tau_2+\delta)$. Moreover, $g_2=g_1$ outside
$B(\tau_2+\delta)$ and $g_2$ can be chosen to be close to
  the metric $g_1$ in the $C^0-$topology. In particular, if
$\delta>0$ is small enough, then the first eigenvalue of the
Laplacian operator of  $(\mathbf{D},g_2)$ is bounded below by a
constant $C_1>0$ independent of $\delta$. Since $0\ge R_2+6\ge -
2\delta$, by \cite{FS}, if $\delta$ is small enough then  there is
a positive solution $v$ of
$$\Delta_{g_2}v-\frac18(R_2+6)v=0.$$

We want to conformally deform $g_2$ to an AH metric  with constant
scalar curvature -6. To do this, for any $k>0$, consider the
following boundary value problem
\begin{equation}\label{3rdmetric}
\left\{
\begin{array}{ll}
 \Delta_{g_2} u_k -\frac18 R_2u_k -\frac34 u_k ^5 =0  , & \hbox{in $\hat B(k)$;} \\
    u_k |_{\partial \hat B (k)}=1, \\
\end{array}
\right.
\end{equation}
where $\hat B(k)$ is the geodesic ball with respect to
$g_2$ of radius $k$ with center at the origin of
$\mathbf{D}$. By rescaling $v$ in $\hat B(k)$ we may assume
that $v>1$ in $\hat B(k)$. Then we have $$\Delta_{g_2} v
-\frac18 R_2v -\frac34 v^5 \le 0$$ in $\hat B(k)$ and the
constant function $v_1=1$ satisfies
$$
\Delta_{g_2} v_1 -\frac18 R_2v_1 -\frac34 v_1^5 \ge 0.
$$  Here we have used the fact that $R_2\le
-6$. By \cite{Sa} as in the proof of Theorem \ref{AdsSch3},
(\ref{3rdmetric}) has a solution $u_k\ge 1$. Suppose $u_k$ attains
maximum at a point $x_0\in \hat B(k)$, then we have $$\frac18
R_2u_k +\frac34 u_k ^5=\Delta_{g_2}u_k\le0
$$
at $x_0$. Hence $$\max_{\hat B(k)}u_k\le
\max_{\mathbf{D}}(-\frac16 R_2)^\frac14.$$

In particular, $u_k$ are uniformly bounded. By taking a
subsequence if it is necessary, we see that there is a smooth
function $u\ge 1$ on $\mathbf{D}$ satisfying:
$$\Delta_{g_2} u -\frac18 R_2u -\frac34 u^5 =0.$$

We claim that $\lim_{x\to\infty}u(x)=1$. Let $g=u^4g_2$, then the
scalar curvature of $g$ is $-6$. Moreover, outside
$B(\tau_2+\delta)$, $g=u^4g_1=u^4\phi_\e^4ds^2_{\H^3}$. $u\phi_\e$
satisfies:
$$
\Delta_{\H^3}(u\phi_\e)+\frac34 u\phi_\e[1-(u\phi_\e)^4]=0.
$$
Use the functions in defining the anti-de
Sitter-Schwarzschild metric in Proposition \ref{AdsSch1} as
comparison functions we conclude that $u\phi_\e\to 1$ as
$x\to\infty$. This proves the claim.

We want to prove that $g$ is an AH metrics with mass $M_g$ such
that $|M_g-M|=C(\delta)$ with $C(\delta)\to0$ as $\delta\to0$ and
that $g$ has surfaces $S_1, S_2, S_3$ as in the theorem. Note that
$$
|u-1|\le C(\delta)
$$
on $B(\tau_2+\delta)$ as in the proof of Theorem \ref{AdsSch2}. Here
$C(\delta)\to0$ as $\delta\to0$. By this and Lemma 1.3, and together
with the standard theory of elliptic partial differential equations,
we see that
\begin{equation}
\|u\phi_\e -\phi_\e\|_{C^{2, \alpha}(\Omega)} \leq C( \Omega,
\delta).
\end{equation}
Here $C( \Omega, \delta)$ is a positive constant which depend only
on $\Omega$ and $\delta$, and $C(\Omega, \delta)\to0$ as
$\delta\to0$. The results will be consequences of Lemma
\ref{maximum} and the following Lemmas \ref{AH}, \ref{mass1}, and
\ref{cmc1}.
\end{proof}

Before we state and prove the lemmas, let us consider the
metric:
\begin{equation}\label{AH0}
\begin{split}
g&=
u^4ds^2_{\H^3}\\
&=\frac{4u^4}{(1-|x|^2)^2}(dr^2+r^2d\sigma^2)\\
&=\frac{4u^4}{(1-|x|^2)^2}e^{2t}(dt^2+
d\sigma^2)\\
&=\frac{4u^4}{(1-|x|^2)^2}g_0
\end{split}
\end{equation}
where $r=|x|$, $r=e^t$, and
$d\sigma^2=h_{\alpha\beta}d\sigma_\alpha d\sigma_\beta$ is the
standard metric on $\SS^2$ and $g_0$ is the Euclidean metric.
Suppose $g$ has constant scalar curvature $-6$ and $u\to1$ as
$|x|\to1$,  then by \cite{ACF}, $u$ is smooth as a function of $x$
up to $|x|=1$.

\begin{lemm}\label{AH} Assume that $|u-1|\le Ce^{-3d_{\H^3}(x,0)}$.
Then $g$ is AH.  \end{lemm}
\begin{proof} Let $\rho=\frac{1-r^2}{2u^2}$. Then $\rho$ is
smooth up to $|x|=1$ and   $|\nabla_0\rho|=1$ at $|x|=1$, where
$\nabla_0$ is the Euclidean gradient. As in \cite [Lemma 5.3]{AD}
(see also p.102 in \cite{CoH}), let $\theta$ be the solution of the
equation, with $\theta=1$ at $r=1$:
\begin{equation}\label{AH1}
   \rho|\nabla_0 \theta|^2+2\theta\langle
   \nabla_0\theta,\nabla_0\rho)=\theta^4\rho+\theta^2a
\end{equation}
with $\theta=1$ at $t=0$, where $a\rho=1-|\nabla_0\rho|^2$
is a smooth function, and $\nabla_0$ is with respect to the
Euclidean metric $g_0$.   Let $f$ be such that $\sinh
f=\theta\rho$. Then
$$
g=\sinh^{-2} f (df^2+g{_{_f}})
$$
where $df^2+g{_{_f}}=\theta^2 g_0=\hat g$, and $g{_{_f}}$ is the
restriction of $\hat g$   on the level surface $f$=constant. Near
$r=1$, i.e., $t=0$, the level surface is a graph of a function
$t=t(\sigma)$ where $\sigma\in\SS^2$. Since
$f(t(\sigma_1,\sigma_2),\sigma_1,\sigma_2)=c$, we have
\begin{equation}\label{AH2}
  f_tt_{\sigma_\alpha}+f_{\sigma_\alpha}=0.
\end{equation}
Since $|\hat\nabla f|^2=1$, here $\hat\nabla$ is with respect to the
metric $\hat g$,  we have
\begin{equation}\label{AH4}
    \theta^{-2}e^{-2t}(f_t^2+h^{\alpha\beta}f_{\sigma_\alpha}f_{\sigma_\beta})=1.
\end{equation}
In local coordinates $(\sigma_1, \sigma_2)$, the metric
$\gamma=g_f$ on the level surface is:
\begin{equation}\label{AH3}
\begin{split}
    \gamma_{_{\alpha\beta}}&=\gamma(\frac{\p}{\p\sigma_\alpha},
    \frac{\p}{\p\sigma_\beta})\\
   &=\theta^2e^{2t}\lf(t_{\sigma_\alpha}t_{\sigma_\beta}
   +h_{\alpha\beta}\ri)\\
   &=\theta^2e^{2t}\lf(\frac{f_{\sigma_\alpha}f_{\sigma_\beta}}{f_t^2}
   +h_{\alpha\beta}\ri)\\
   &=\frac{f_{\sigma_\alpha}f_{\sigma_\beta}}
   {1-h^{\xi\zeta}f_{\sigma_\xi}f_{\sigma_\zeta}}+
   (\theta^2e^{2t}-1)h_{\alpha\beta}+h_{\alpha\beta}.
   \end{split}
\end{equation}
We want to prove that the last term is $h_{\alpha\beta}+O(t^3)$
near $t=0$.

Since $|u-1|=O(e^{-3d(x,0)})$, $|u-1|=O(t^3)$, at $t=0$, and $u$
is smooth up to the boundary,   we have $u_t=u_{tt}=0$ at $t=0$.
Hence at $t=0$,
\begin{equation}\label{AH4}
\rho_t=-1, \rho_{tt}=-2, \rho_{ttt}=-4, \rho_{tttt}=8(-1+u_{ttt}),
\rho_{t_{\sigma_\alpha}}= \rho_{tt_{\sigma_\alpha}}=
\rho_{ttt_{\sigma_\alpha}}=0
\end{equation}

Now
$$
|\nabla_0\rho|^2=e^{-2t}(\rho_t^2+h^{\alpha\beta}
\rho_{\sigma_\alpha}\rho_{\sigma_\beta})=e^{-2t}A$$ we have at
$t=0$,
$$A=1,\ A_t=4, \ A_{tt}=16, A_{ttt}=64-16u_{ttt}
$$
Hence at $t=0$,
\begin{equation}\label{AH5}
   \begin{split}(1-|\nabla_0\rho|^2)_t&=-e^{-2t}[-2A+A_t]=-2,\\
   (1-|\nabla_0\rho|^2)_{tt}&=-e^{-2t}[4A-4A_t+A_{tt}]=-4\\
   (1-|\nabla_0\rho|^2)_{ttt}&=-e^{-2t}[-8A+12A_t-6A_{tt}+A_{ttt}]=8-16u_{ttt}.
   \end{split}
\end{equation}
Now $a\rho=1-|\nabla_0\rho|^2$, at $t=0$, we have
\begin{equation}\label{AH6}
    a=2,\ a_t=0,\ a_{tt}=-\frac{16}3u_{ttt}
\end{equation}
By (\ref{AH1}),
\begin{equation}
   \rho(\theta_t^2+h^{\alpha\beta}\theta_{\sigma_\alpha}\theta_{\sigma_\beta})
   +2\theta(\theta_t\rho_t+h^{\alpha\beta}\theta_{\sigma_\alpha}\rho_{\sigma_\beta})=
   e^{2t}(\theta^4\rho+\theta^2a)
\end{equation}
Note that at $t=0$, $\theta=1$, $\theta_{\sigma_\alpha}=0$. Hence
at $t=0$, $\theta_t=-1$. Now
\begin{equation}
   \begin{split}
   \rho_t&(\theta_t^2+h^{\alpha\beta}\theta_{\sigma_\alpha}
   \theta_{\sigma_\beta})+\rho(\theta_t^2+h^{\alpha\beta}\theta_{\sigma_\alpha}
   \theta_{\sigma_\beta})_t
   \\&+2\theta_t(\theta_t\rho_t+h^{\alpha\beta}
   \theta_{\sigma_\alpha}\rho_{\sigma_\beta})
   +2\theta(\theta_{tt}\rho_t+\theta_t\rho_{tt}+[h^{\alpha\beta}
   \theta_{\sigma_\alpha}\rho_{\sigma_\beta}]_t)\\
   &=2e^{2t}(\theta^4\rho+\theta^2a)+e^{2t}[(\theta^4)_t\rho+\theta^4\rho_t+
   2\theta\theta_ta+\theta a_t]
   \end{split}
\end{equation}
So at $t=0$, $\theta_{tt}=1$. Here we have used the fact that
$\theta_{t\sigma_\alpha}=\rho_{t\sigma_\alpha}=0$ at $t=0$.
\begin{equation}
   \begin{split}
   \rho_{tt}&(\theta_t^2+h^{\alpha\beta}\theta_{\sigma_\alpha}
   \theta_{\sigma_\beta})+2\rho_t(2\theta_t\theta_{tt}+[h^{\alpha\beta}\theta_{\sigma_\alpha}
   \theta_{\sigma_\beta}]_t)+\rho(\theta_t^2+h^{\alpha\beta}\theta_{\sigma_\alpha}
   \theta_{\sigma_\beta})_{tt}\\
   &+2\theta_{tt}(\theta_t\rho_t+h^{\alpha\beta}\theta_{\sigma_\alpha}\rho_{\sigma_\beta})+4\theta_t(\theta_{tt}\rho_t+\theta_{t}\rho_{tt}+[h^{\alpha\beta}
   \theta_{\sigma_\alpha}\rho_{\sigma_\beta}]_t)\\
   &+2\theta(\theta_{ttt}\rho_t+2\theta_{tt}\rho_{tt}+\theta_{t}\rho_{ttt}+[h^{\alpha\beta}
   \theta_{\sigma_\alpha}\rho_{\sigma_\beta}]_t)\\
   &= e^{2t}\{4( \theta^4 \rho+\theta^2a)+ 4[(\theta^4)_t\rho
   + \theta^4 \rho_t+2\theta\theta_t a+\theta^2a_t]\\
   &+(\theta^4)_{tt}\rho
   +8\theta^3\theta_t\rho_t+\theta^4\rho_{tt}+(2\theta_t^2
   +2\theta\theta_{tt})a+2(\theta^2)_ta_t+\theta^2a_{tt}\}
   \end{split}
\end{equation}
Hence at $t=0$, $\theta_{ttt}=-1+\frac18u_{ttt}$.
\begin{equation}\label{AH6}
   \begin{split}
   (\theta^2e^{2t}-1)_t&=2\theta\theta_te^{2t}+2e^{2t}\theta^2=0\\
(\theta^2e^{2t}-1)_{tt}&=e^{2t}
\lf(2\theta\theta_{tt}+2\theta_t^2+4\theta^2+8\theta\theta_t\ri)=0\\
 (\theta^2e^{2t}-1)_{ttt}&=
 2e^{2t}\lf(2\theta\theta_{tt}+2\theta_t^2+4\theta^2+8\theta\theta_t\ri)\\
 &\quad+
 e^{2t}\lf(2\theta\theta_{ttt}+6\theta_t\theta_{tt}+
 8\theta\theta_t+8\theta_t^2+8\theta\theta_{tt}\ri)=2\theta_{ttt}+2
   \end{split}
\end{equation}
Since $\sinh f=\theta\rho$, $$ \cosh f
f_t=\theta_t\rho+\theta\rho_t=-1 $$ at $t=0$. Hence $f_t=-1$ at
$t=0$ and so $f_{\sigma_\alpha t}=0$ at $t=0$.
$$
\cosh f f_{tt}+\sinh f
f_t^2=\theta_{tt}\rho+2\theta_t\rho_t+\theta\rho_{tt}=-4.
$$
So $f_{tt}=-4$ at $t=0$ and   $f_{tt\sigma_\alpha}=0$ at $t=0$.
From these and (\ref{AH3}), we conclude that
$$
\gamma_{\alpha\beta}=h_{\alpha\beta}+O(t^3),
$$
near $t=0$. Hence $g$ is AH because $u$ is smooth at $|x|=1$.
\end{proof}
To get an expression for the mass, let us compute $f_{ttt}$, we
have
\begin{equation}\label{AH7}
\begin{split}
\cosh f f_{ttt}+3\sinh f f_tf_{tt}+ \cosh f
f_t^3&=\theta_{ttt}\rho+3\theta_{tt}\rho_t+
3\theta_{t}\rho_{tt}+\theta\rho_{ttt}.
\\
&=-3+6-4 \\
&=-1.
\end{split}
\end{equation}
So $f_{ttt}=0$ at $t=0$. Hence if
$$\gamma_{\alpha\beta}=h_{\alpha\beta}+\frac{f^3}{3}\eta_{\alpha\beta}+O(f^4)
$$
then
$$
\gamma_{\alpha\beta}=h_{\alpha\beta}-\frac{t^3}{3}\eta_{\alpha\beta}+O(t^4).
$$
Hence
\begin{equation}\label{AH8}
\begin{split}
\eta_{\alpha\beta}&=-\frac12(\gamma_{\alpha\beta})_{ttt}\\
&=-\frac12(\theta^2e^{2t}-1)_{ttt}h_{\alpha\beta}\\
&=-(\theta_{ttt}+1)h_{\alpha\beta}\\
&=-\frac18u_{ttt}h_{\alpha\beta}
\end{split}
\end{equation}
evaluated at $t=0$.
\begin{lemm}\label{mass1}
$\text{tr}_{h}\gamma=-\frac14u_{ttt}$. Hence suppose
$g_1=u^4_1ds^2_{\H^3}$ and $g_2=u_2^4ds^2_{\H^3}$ be two
metrics defined outside some compact set of $\H^3$ such
that $u_1$ and $u_2$ are smooth up to $\p B(1)$ if we use
ball model for $\H^3$. Moreover,  assume that
$|u_1(x)-1|+|u_2(x)-1|\le Ce^{-3d(x,0)}$. Suppose $u_1$ and
$u_2$ are close in the sense that $|u_1-u_2|\le \epsilon
e^{-3d(x,0)}$. Then there is an absolute constant $C_1$
such that $|M_1-M_2|\le C_1\epsilon$, here $M_1$, $M_2$ is
the mass of $g_1$, $g_2$ respectively.
\end{lemm}
\begin{proof} To prove the second part, with the same notation  as in Lemma
\ref{AH}, we have $|(u_1)_{ttt}-(u_2)_{ttt}|\le C_2\epsilon$ at
$t=0$ for some absolute constant $C_2$. The result follows from
first part and the definition of mass.
\end{proof}
Since $g_2$ in the proof  Theorem \ref{AdsSch4} is anti-de
Sitter-Schwarzchild outside $B(\tau_2+\delta)$, we may assume
$g_2=\phi^4(dr^2+r^2d\sigma^2)$ for some $r\geq \delta>0$, where
$d\sigma^2=h_{\alpha\beta}d\sigma_\alpha d\sigma_\beta$ is the
standard metric on $\SS^2$. Without loss of generality, we may
assume $(\mathbf{D}_{1-\delta}\setminus \mathbf{D}_\delta , g_0)$
containing the compact surface with mean curvature $\pm2$ and $0$,
here and in the sequel, the mean curvature is always with respect
to the outward unit normal vector, and for simplicity, we denote
$\mathbf{D}_{1-\delta}\setminus \ol{\mathbf{D}}_\delta$ by $N$

Obviously, it is enough to show

\begin{lemm}\label{cmc1} Let $g_2= \phi^4 (dr^2+r^2d\sigma^2)$ be as
above,
which is a Riemannian metric on $N$.
 Then there is an $\e>0$ such that for any $\tilde\phi$ with
$\|\phi-\tilde\phi\|_{C^{2, \alpha}(N)}\leq \e$, then there are
compact surfaces in $(N,g)$ with   mean curvature equal to $\pm2$
and $0$, here $g=\tilde\phi^4 (dr^2+r^2d\sigma^2)$.
\end{lemm}

\begin{proof} We will use implicit function theorem. Let us discuss the case that $H=0$ and $H=\pm2$ together,
we adopt the coordinates $(r,\sigma)$ on $N$, here $\sigma \in
\SS^2$. Consider the Banach spaces $\frak B_1=C^{2, \alpha}(\ol
N)$, $\frak B_2=C^{2, \alpha}(\mathbf{\SS}^2)$ and $\frak
B_3=C^{\alpha}(\SS^2) $ and let $\frak U$ be the open set in
$\frak B_1\times \frak B_2$ consisting of $(\Phi, v)$ such that
$\Phi>0$ and $\delta <v<1-\delta$. For $(\Phi, v)\in \frak U$, let
define $H(\Phi,v)$ to be the mean curvature of the surface given
by $(\sigma,v(\sigma))$ in $(N, \Phi^4(dr^2+\sigma^2d\sigma^2)$.
Then $H:\frak U\to \frak B_3$.

 We want to compute the differential at  $\Phi=\phi=\phi(r)$ and
 $v=c$=constant with $\delta<c<1-\delta$. Let $\tilde
\nabla=\nabla_{\SS^2}$. Consider a surface  given by
$(\sigma,v(\sigma))$ and let $f(r,\sigma)=v(\sigma)-r$. Then the
surface is given by the level surface  $f=0$. Then
\begin{equation}\label{meancurvature1}
   \left\{%
\begin{array}{ll}
     \nabla f& = \phi^{-4}(-\frac{\p}{\p r}+r^{-2} \nabla v); \\
    \Delta f&=-\frac{1}{\phi^6r^2}(\phi^2r^2)_r+\phi^{-4}r^{-2}\tilde\Delta v;  \\
    |\nabla f|^2&=\phi^{-4}(1+r^{-2}|\tilde\nabla v|^2)=\phi^{-4}\psi^2 \\
\end{array}%
\right.
\end{equation}
where $\nabla$, $\Delta$ are the gradient and Laplacian with respect
to $g_2$, $\tilde \Delta$ is the Laplacian on $\SS^2$ and
$\psi=(1+r^{-2}|\tilde\nabla v|^2)^\frac12$. Therefore the mean
curvature of the level surface $f=0$ is
\begin{equation}\label{meancurvature2}
   \begin{split}
    H(\phi, v)&=div\lf(\frac{\nabla f}{|\nabla f|}\ri)\\
    &=\frac{\Delta f}{|\nabla f|}+\langle \nabla (\frac1{|\nabla
    f|}),\nabla f \rangle\\
    &=-\phi^{-4}\psi^{-1}r^{-2}(\phi^2
    r^2)_r+\phi^{-2}\psi^{-1}r^{-2}\tilde \Delta
    v-\phi^{-4}(\phi^2\psi^{-1})_r\\
    &\ \ +\phi^{-2}r^{-2}\langle
    \tilde\nabla(\psi^{-1}),\tilde \nabla v\rangle\\
    &=-\phi^{-4}\lf[2\psi^{-1}(\phi^2)_r+2\psi^{-1}r^{-1}\phi^2+\phi^2(\psi^{-1})_r\ri]
    \\
    &\ \ +\phi^2\psi^{-1}r^{-2}\tilde\Delta v+\phi^{-2}r^{-2}\langle
    \tilde\nabla(\psi^{-1}),\tilde \nabla v\rangle
    \end{split}
\end{equation}
evaluated at $r=v$. By the expression of $H$, it is easy to see that
$H$ is $C^1$

We want to compute $\frac{d}{d t}H(\phi, c+t\eta)|_{t=0}$, where
$c$ is a constant. It is easy to see that
\begin{equation}\label{meancurvature2}
   \frac{d}{dt}\psi(c+t\eta,\sigma)|_{t=0}= \frac{d}{dt}\frac{\p}{\p
   r}(\psi(c+t\eta,\sigma))|_{t=0}=0.
\end{equation}
Since $c$ is a constant, we have
\begin{equation}\label{meancurvature3}
\frac{d}{dt}\lf[\phi^2\psi^{-1}r^{-2}\tilde\Delta
v+\phi^{-2}r^{-2}\langle
    \tilde\nabla(\psi^{-1}),\tilde \nabla
    v\rangle\ri]|_{t=0}=\phi^{-2}c^{-2}\tilde\Delta\eta
\end{equation}
Hence
\begin{equation}\label{meancurvature4}
\begin{split}
\frac{d}{d
t}H(\phi, c+t\eta)|_{t=0}&=4\phi^{-5}\phi_r\lf[4\phi\phi_r+2c^{-1}\phi^2\ri]\eta\\
&\ \ -
\phi^{-4}\lf[4\phi_r^2+4\phi\phi_{rr}-2c^{-2}\phi^2+4c^{-1}\phi\phi_r\ri]\eta\\
&\ \ +\phi^{-2}c^{-2}\tilde\Delta\eta.
\end{split}
\end{equation}
Since the metric $g_2$ has constant scalar curvature -6 on $N$ and
the surface $r=c$ has constant mean curvature $H$,
$$
H=\phi^{-2}(\frac 2r+4\phi^{-1}\phi_r)=\phi^{-2}(
2c^{-1}+4\phi^{-1}\phi_r) $$ and so
$$
\phi_r=\frac14\phi^3 H-\frac12\phi c^{-1}.
$$
Also
$$
8\Delta_0\phi=6\phi^5,
$$
where $\Delta_0$ is the Euclidean Laplacian. Hence
$$
\phi_{rr}=-\frac 2r\phi_r+\frac34\phi^5=-2c^{-1}\lf(\frac14\phi^3
H-\frac12\phi c^{-1}\ri)+\frac34\phi^5.
$$
Therefore
\begin{equation}\label{meancurvature4}
  \begin{split}
    \frac{d}{d
t}&H(\phi,c+t\eta)|_{t=0}\\
&=
\eta\phi^{-4}\lf[4\phi^3H\phi_r-4\phi_r^2-
4\phi\lf(-2c^{-1}\phi_r+\frac34\phi^5\ri)+2c^{-2}\phi^2-4c^{-1}\phi\phi_r\ri]\\
 &  \ \ +\phi^{-2}c^{-2}\tilde\Delta\eta. \\
  &= \eta\phi^{-4}\lf[4\phi^3H\phi_r-4\phi_r^2+
4\phi c^{-1}\phi_r+2c^{-2}\phi^2-3\phi^6 \ri]\\
 &  \ \ +\phi^{-2}c^{-2}\tilde\Delta\eta\\
 &=\eta\phi^{-4}\lf[4\phi^3H\phi_r-4\phi_r^2+
4\phi c^{-1}\phi_r+2c^{-2}\phi^2-3\phi^6 \ri]\\
 &  \ \ +\phi^{-2}c^{-2}\tilde\Delta\eta\\
 &=\eta\phi^{-4}\lf[4\phi^3H\phi_r-4\phi_r^2+\phi^4c^{-1}H-3\phi^6 \ri]\\
 &  \ \ +\phi^{-2}c^{-2}\tilde\Delta\eta\\
 &=\eta\phi^{-4}\lf[ \phi^6H^2 -2\phi^4c^{-1}H-\frac14\phi^6H^2+\phi^4c^{-1}H-\phi^2c^{-2}
 +\phi^4c^{-1}H-3\phi^6 \ri]\\
 &  \ \ +\phi^{-2}c^{-2}\tilde\Delta\eta\\
 &=\eta\phi^{-4}\lf[ \frac34\phi^6H^2 -\phi^2c^{-2}-3\phi^6 \ri]\\
 &  \ \ +\phi^{-2}c^{-2}\tilde\Delta\eta\\
 &=\phi^{-2}c^{-2}\tilde\Delta\eta- \Theta\eta
\end{split}
\end{equation}
where $\Theta\ge \phi^{-2}c^{-2} $ is a positive function if
$|H|\le 2$.

Let $c$ be such that $H(\phi, c)=2$, then $\delta <c<1-\delta$.
$\frac{d}{d t}H(\phi,c)|_{t=0}$: $C^{2,
\alpha}(\mathbf{\SS}^2)\rightarrow C^{\alpha}(\SS^2)$ is
bijective, thus, by implicit function theorem, we see that there
is $\e>0$ so that for any $\|\phi-\tilde\phi\|_{C^{2,
\alpha}(N)}\leq \e$, there is a smooth function $v$ on $\SS^2$
with $\|v-c\|_{C^{2, \alpha}(\SS^2)}\leq \e$ such that
$H(\tilde\phi, v)=2$. Thus, there is a compact surface in $(N,g)$
with the mean curvature equal to $2$. by the same arguments, one
may show there are compact surfaces with mean curvature $-2$ and
$0$.
\end{proof}
Together with (4.2), we see that there exists surfaces $S_1$, $S_2$,
$S_3$ in the manifolds that we constructed. Note that the surfaces
are diffeomorphic to $\SS^2$ and are close to the constant mean
curvature surfaces in the metric $g_2$. Thus, we finish to prove
Theorem 3.1.

\bibliographystyle{amsplain}

\bibliographystyle{amsplain}

\end{document}